\documentclass[12pt,reqno]{amsart}
\usepackage{fullpage}
\usepackage[top=1.5in, bottom=2in, left=2.5in, right=2.5in]{geometry}
\usepackage{amssymb,amsfonts,mathtools}
\usepackage[all,arc]{xy}
\usepackage{enumerate}
\usepackage{mathrsfs}
\usepackage{amsmath}
\usepackage{amsthm}
\usepackage{tikz}
\usepackage{tikz-cd}
\usepackage[utf8]{inputenc}
\usepackage{lipsum}
\usepackage{extpfeil}
\usepackage{colonequals}
\usepackage{graphicx}
\usepackage{subfig}
\usepackage{enumitem}
\usepackage{array}
\usepackage{xcolor}
\colorlet{RED}{red}
\usepackage{bm, bbm}
\usepackage{url}

\newtheorem{thm}{Theorem}[section]

\newtheorem*{conj*}{Conjecture}

\newtheorem{prop}[thm]{Proposition}
\newtheorem{lem}[thm]{Lemma}

\theoremstyle{definition}

\theoremstyle{remark}

\DeclareMathOperator{\ZZ}{\mathbb{Z}}

\numberwithin{equation}{section}

\newcommand{\lp}{\left(}
\newcommand{\rp}{\right)}

\newcommand{\ou}{\mathrm{ou}}
\newcommand{\eu}{\mathrm{eu}}
\newcommand{\od}{\mathrm{od}}
\newcommand{\ed}{\mathrm{ed}}
\newcommand{\x}{\mathrm{x}}
\newcommand{\y}{\mathrm{y}}
\newcommand{\z}{\mathrm{z}}
\newcommand{\w}{\mathrm{w}}
\renewcommand{\o}{\mathrm{o}}
\newcommand{\e}{\mathrm{e}}
\renewcommand{\d}{\mathrm{d}}
\renewcommand{\u}{\mathrm{u}}

\renewcommand{\pmod}[1]{\ \lp  \mathrm{mod} \, #1 \rp }

\allowdisplaybreaks

\author{Kathrin Bringmann}
\address{University of Cologne, Department of Mathematics and Computer Science, Weyertal 86-90, 50931 Cologne, Germany}
\email{kbringma@uni-koeln.de}

\author{Catherine Cossaboom}
\address{University of Virginia, Charlottesville, VA 22903, USA}
\email{qkb9us@virginia.edu}

\author{William Craig}
\address{Math Department, United States Naval Academy}
\email{wcraig@usna.edu}

\usepackage{tikz-cd}
\usepackage{float}
\usepackage{array}
\usepackage{comment}
\subjclass[2020]{}
\keywords{Partitions, parity, overpartitions, $q$-hypergeometric series.}

\title[Overpartitions with parts separated by parity]{Overpartitions with parts separated by parity}

\begin{document}
	
	\begin{abstract}
		In this paper, we generalize Andrews' partitions separated by parity to overpartitions in two ways. We investigate the generating functions for $16$ overpartition families whose parts are separated by parity, and we prove various $q$-series identities for these functions. These identities include relations to modular forms, $q$-hypergeometric series, and mock modular forms.
	\end{abstract}
	
	\maketitle
	
	\section{Introduction and Statement of Results}
	
	We recall that a {\it partition} of an integer $n$ is a non-increasing sequence $\lambda = \lambda_1+ \dots+ \lambda_\ell$ of positive integers whose parts, denoted by $\lambda_j$ ($1 \leq j \leq \ell$), sum up to $n$. Restrictions on partitions concerned with the parity of the parts in the partition have played a major role in the history of partition theory; see e.g. Andrews' article \cite{AndrewsParity} for an excellent overview of many ways parity has appeared in partition theory. In Andrews' recent works in this direction \cite{Andrews2018,Andrews2019}, he introduced the notion of ``partitions with parts separated by parity'' in connection with mock theta functions. A partition $\lambda$ has {\it parts separated by parity}\footnote{Partitions into only odd parts or only into even parts are permitted.} if each even part $\lambda_j$ is larger than each odd part $\lambda_k$, or if each odd part $\lambda_j$ is larger than each even part $\lambda_k$. Andrews used the notation $p_{\x\y}^{\z\w}(n)$ to count various kinds of partitions of $n$ with parts separated by parity. In his notation, we let $\{\x,\z\}=\{\e,\o\}$, where the case $\z = \e, \x = \o$ (resp. $\z = \o, \x = \e$) signifies that even parts are larger than odd parts (resp. that odd parts are larger than even parts). Andrews then allowed $\y$ and $\w$ to denote either $\u$ or $\d$, which represent unrestricted and distinct, respectively, and signify whether parts of parity $\x$ and $\z$ are unrestricted or must be distinct, respectively. Andrews studied the eight possible functions built this way; e.g. $p_{\od}^{\eu}(n)$ counts the number of partitions of $n$ with even parts larger than odd parts which have all odd parts distinct.
	
	Andrews' papers \cite{Andrews2018,Andrews2019} and the follow-up work of the first author and Jennings-Shaffer \cite{BringmannJenningsShaffer} focused on the generating functions
	\begin{align*}
		F_{\x\y}^{\z\w}(q) := \sum_{n \geq 0} p_{\x\y}^{\z\w}(n) q^n.
	\end{align*}
	These generating functions have connections to a wide variety of modular-type objects. To state some of these results, we define the \emph{$q$-Pochhammer symbol}
	\begin{align*}
		\lp a \rp_n = \lp a;q \rp_n := \prod_{k=0}^{n-1} \lp 1 - aq^k \rp, \ \ \ \ \ \lp a_1, a_2, \dots, a_m; q \rp_n := \prod_{k=1}^m \lp a_k; q \rp_n,
	\end{align*} for $n \in \mathbb{N}_0 \cup \{\infty\},\,m\in\mathbb{N}$. As a first example, Andrews proved \cite[equation (2.1)]{Andrews2018} that
	\begin{align*}
		F_{\eu}^{\ou}(q) = \dfrac{1}{(1-q)\lp q^2;q^2 \rp_\infty},
	\end{align*}
	which is essentially a modular form of weight $-\frac{1}{2}$. The papers \cite{Andrews2019,BringmannJenningsShaffer} both proved (see also notation and comments from \cite{BringmannCraigNazaroglu1,BringmannCraigNazaroglu2}) that
	\begin{align*}
		F_{\od}^{\eu}(q) = \dfrac{1}{\lp q^2;q^2 \rp_\infty}\lp 1 - \dfrac{\sigma(-q)}{2} + \dfrac{\lp -q;-q \rp_\infty}{2} \rp,
	\end{align*}
	where
	\begin{align*}
		\sigma(q) := \sum_{n \geq 0} \dfrac{q^{\frac{n(n+1)}{2}}}{\lp -q;q \rp_n}
	\end{align*}
	is the famous Ramanujan $\sigma$-function from his lost notebook \cite{Ra}. It was shown by Andrews--Dyson--Hickerson \cite{AndrewsDysonHickerson} and Cohen \cite{Cohen} that $\sigma(q)$ is related to Maass forms. The variety of generating functions in \cite{Andrews2019} also connects to modular forms and mock modular forms. Andrews' work has given rise to many other examples, looking at further generating function identities \cite{BringmannJenningsShaffer}, asymptotic formulas \cite{BringmannCraigNazaroglu1,BringmannCraigNazaroglu2}, and connections with mock theta functions \cite{FuTang}. Although the definitions of each $p_{\x\y}^{\z\w}(n)$ are quite similar, the properties that arise from their generating functions are very different (including modular forms, mock modular forms, false modular forms, and mock Maass theta functions like $\sigma(q)$); this suggests a rich theory yet to be fully understood. 
    
	In this paper, we extend this framework to overpartitions. An {\it overpartition} \cite{CorteelLovejoy} is a partition where the first instance of each part may possibly be overlined. For instance, there are 8 overpartitions of 3, given by
	\begin{align*}
		3,\ \overline 3,\ 2+1,\ \overline 2 + 1,\ 2 + \overline 1,\ \overline 2 + \overline 1,\ 1 + 1 + 1,\ \overline 1 + 1 + 1.
	\end{align*}
	Overpartitions have played numerous roles in $q$-series and combinatorics (see e.g.\! \cite{BP,Lo} and the references therein), mathematical physics (see e.g.\! \cite{FJM}), symmetric functions (see e.g.\! \cite{Br}), representation theory (see e.g.\! \cite{KK}), and algebraic number theory (see e.g.\! \cite{Lo2}). The goal of this paper is to pursue this theme in the context of partitions with parts separated by parity.
	
	We define two variations of overpartitions with parts separated by parity. In both cases, we say that an overpartition $\lambda$ has {\it parts separated by parity}\footnote{As in the standard case, overpartitions into only odd parts or only even parts are permitted.} if each even $\lambda_j$ is larger than each odd $\lambda_k$, or if each odd $\lambda_j$ is larger than each even $\lambda_k$. 
	Because allowing the distinct parts to be overlined introduces powers of two that are both complicated and unenlightening, we introduce an additional constraint in order to create an alternative family of overpartitions. If we impose the restriction that any parts required to be distinct cannot be overlined, then we call these {\it modified}. We let $\overline{p}_{\x\y}^{\z\w}\lp n\rp$ and $\underline{p}_{\x\y}^{\z\w}\lp n\rp$ denote the functions which count the number of overpartitions with parts separated by parity and the modified variation, respectively. For instance, we have that $\overline{p}_{\od}^{\eu}(3) = 6$ and $\underline{p}_{\od}^{\eu}(3) = 3$, respectively, since the corresponding overpartitions are
	\begin{align*}
		3 ,\ 2+1,\ \overline 3,\ \overline 2 + 1,\ 2 + \overline 1,\ \overline 2 + \overline 1
	\end{align*}
	and
	\begin{align*}
		3,\ 2+1,\ \overline 2 + 1.
	\end{align*}
	The difference between the two cases is that the overpartitions counted by $\overline{p}_{\od}^{\eu}(3)$ are permitted to have overlines on both the even and odd parts, whereas the overpartitions counted by $\underline{p}_{\od}^{\eu}(3)$ are only permitted to have overlines on the even parts.
	
	In this paper, we consider the 16 generating functions
	\begin{align*}
		\overline{F}_{\x\y}^{\z\w}(q)  := \sum_{n \geq 0} \overline{p}_{\x\y}^{\z\w}\lp n\rp  q^n, \ \ \ \ \ \underline{F}_{\x\y}^{\z\w}(q) := \sum_{n \geq 0} \underline{p}_{\x\y}^{\z\w}\lp n\rp  q^n.
	\end{align*}
	We also need the {\it Jacobi theta function}
	\begin{align*}
		\Theta\lp \tau\rp  := \sum_{n \in \ZZ} q^{n^2} \ \ \ \ \ (q:=e^{2\pi i\tau}).
	\end{align*}
	This is a modular form of weight $\frac{1}{2}$.
	Our first results are formulas for the generating functions $\overline{F}_{\x\y}^{\z\w}(q)$.

	\begin{thm} \label{T: Standard}
		\begin{align*}
			\intertext{{\rm(1)} We have}
			\overline{F}_{\eu}^{\ou}(q) &= \frac{\lp -q^2; q^2\rp _\infty}{2\lp q^2; q^2\rp _\infty} \lp 1 + \Theta^2(\tau) \rp. 
			\intertext{{\rm(2)} We have}
			\overline{F}_{\ou}^{\eu}(q) &= \dfrac{\lp -q^2;q^2 \rp_\infty}{\lp q^2;q^2 \rp_\infty} + \frac{2q}{1-q} \frac{\lp -q^2; q^2\rp _\infty}{\lp q^2; q^2\rp _\infty} \sum_{n \geq 0} \frac{\lp -q, q^2; q^2\rp _{n}}{\left(q^3,-q^2; q^2\right)_{n}} q^{2n}. 
			\intertext{{\rm(3)} We have}
			\overline{F}_{\ed}^{\od}(q) &= \lp -2q;q^2 \rp_\infty - \dfrac{q}{1-q} \lp -2q^2;q^2 \rp_\infty + \dfrac{q}{1-q} \lp -2q;q^2 \rp_\infty. 
			\intertext{{\rm(4)} We have}
			\overline{F}_{\od}^{\ed}(q) &= \lp -2q^2;q^2 \rp_\infty + \dfrac{q}{1-q} \lp 3 \lp -2q^2;q^2 \rp_\infty - \lp -2q;q^2 \rp_\infty \rp. 
			\intertext{{\rm(5)} We have}
			\overline{F}_{\eu}^{\od}(q) &= \dfrac{\lp -q^2;q^2 \rp_\infty}{\lp q^2;q^2 \rp_\infty} \sum_{n \geq 0} \dfrac{2^n q^{n^2}}{\lp -q^2;q^2 \rp_n}. 
			\intertext{{\rm(6)} We have}
			\overline{F}_{\od}^{\eu}(q) &= - \dfrac{\lp -2q, q^2; q^2 \rp_\infty}{\lp -q^2;q^2 \rp_\infty} \sum_{n \geq 0} \dfrac{(-1)^n q^{2n}}{\lp 2q; q^2 \rp_{n+1}} + 2 \sum_{n \geq 0} \dfrac{(-1)^n 2^n q^{n^2+2n}}{\lp 1 + q^{2n+2} \rp \lp 2q; q^2 \rp_{n+1}}.\!\!\!\! 
			\intertext{{\rm(7)} We have }
			\overline{F}_{\ed}^{\ou}(q) &= \frac{\left(-q^2;q^2\right)_\infty}{\left(q^2;q^2\right)_\infty} + 2q\left(-2q^2;q^2\right)_\infty \sum_{n\ge0} \frac{(-1)^n q^{2n}}{\left(2q;q^2\right)_{n+1}} \\ &+ \frac{2q\left(-q;q^2\right)_\infty}{\left(q^3;q^2\right)_\infty} \sum_{n\ge0} \frac{(-1)^n 2^n\left(-q;q^2\right)_n q^{n^2+3n}}{\left(2q,-q^2;q^2\right)_{n+1}}.\\
			\intertext{{\rm(8)} We have}
			\overline{F}_{\ou}^{\ed}(q) &= \lp -2q^2;q^2    \rp_\infty + \sum_{n \geq 0} \dfrac{\lp -q;q^2       \rp_n \lp -2q^{2n+2};q^2 \rp_\infty}{\lp q;q^2       \rp_{n+1}} q^{2n+1}. 
		\end{align*}
	\end{thm}
	
	\noindent{\bf Remark. }\it
	We note a connection to modular forms. In particular, Theorem \ref{T: Standard} {\rm(1)} yields a sum of weakly holomorphic modular forms of weights $-\frac 12$ and $\frac 12$. \rm\\
	
	We next turn to $\underline{F}^{\z\w}_{\x\y}(q)$. We note that in cases where odd and even parts are either both unrestricted or both distinct, the generating functions boil down to either $\overline{F}$ or $F$, respectively; that is,
	\begin{align*}
		\underline{F}_{\eu}^{\ou}(q) = \overline{F}_{\eu}^{\ou}(q), \ \ \ \underline{F}^{\eu}_{\ou}(q) = \overline{F}^{\eu}_{\ou}(q), \ \ \ \underline{F}_{\ed}^{\od}(q) = F_{\ed}^{\od}(q), \ \ \ \underline{F}_{\od}^{\ed}(q) = F_{\od}^{\ed}(q).
	\end{align*}
    These are treated either in Theorem \ref{T: Standard} or in previous works \cite{Andrews2019,BringmannJenningsShaffer}. Thus, we focus on the remaining four examples of interest, in which one parity is restricted and one is unrestricted. In order to state this theorem, we recall the third order mock theta function
    \begin{align*}
		\phi(q) := \sum_{n \geq 0} \dfrac{q^{n^2}}{\lp -q^2;q^2 \rp_n}.
    \end{align*}
    We prove the following result on the remaining modified generating functions.
	
	\begin{thm} \label{T: Modified}
		\begin{align*}
			\intertext{{\rm(1)} We have}
			\,\underline{F}_{\eu}^{\od}(q) &= \dfrac{\lp -q^2;q^2 \rp_\infty}{\lp q^2;q^2 \rp_\infty} \phi(q). 
			\intertext{{\rm(2)} We have}
			\,\underline{F}_{\od}^{\eu}(q) &= \dfrac{\lp -q^2;q^2 \rp_\infty}{\lp q^2;q^2 \rp_\infty} - 2q \lp -q; q^2 \rp_\infty \sum_{n\ge0}\frac{q^n}{\left(-q^2;q^2\right)_{n+1}}\\
			&\hspace{4cm} + 4q \dfrac{\lp -q^2;q^2 \rp_\infty}{\lp q^2;q^2 \rp_\infty} \sum_{n \geq 0} \dfrac{(-1)^n q^{n^2 + 2n}}{\lp 1 + q^{2n+2} \rp \lp q;q^2 \rp_{n+1}}.
			\intertext{{\rm(3)} We have}
			\underline{F}_{\ed}^\ou(q) &= - 2q \lp -q^2; q^2 \rp_\infty \sum_{n \geq 0} \dfrac{q^{n}}{\lp -q^2;q^2 \rp_{n+1}}\\
		&\hspace{4cm} + \dfrac{\lp -q;q^2 \rp_\infty}{\lp q;q^2 \rp_\infty} \sum_{n \geq 0} \dfrac{ (-1)^n \lp -q;q^2 \rp_{n+1} q^{n^2+n}}{\lp -q^2, q; q^2 \rp_{n+1}}. 
            \\
			\intertext{{\rm(4)} We have}
			\,\underline{F}_{\ou}^{\ed}(q) &= \lp -q^2;q^2 \rp_\infty \lp 1 + 2q \sum_{n \geq 0} \dfrac{\lp -q;q^2 \rp_n}{\lp q, -q^2;q^2 \rp_n} q^{2n} \rp.
		\end{align*}
	\end{thm}
    
	The remainder of the paper is structured as follows. In Section \ref{S: Preliminaries}, we provide several known $q$-series identities which we use in the proof of our main results. In Sections \ref{S: Standard} and \ref{S: Modified}, we prove the results regarding standard and modified overpartitions with parts separated by parity, respectively.
	
	\section*{Acknowledgments}
	
	The first and the third author have received funding from the European Research Council (ERC) under the European Union’s Horizon 2020 research and innovation programme (grant agreement No. 101001179). This paper was partially written while the second author was visiting the Max Planck Institute for Mathematics, whose hospitality she acknowledges. She also recognizes support from the Raven Fellowship and the Ingrassia Family Echols Scholars Research Grant. The authors also thank Koustav Banerjee for directing our attention towards Proposition \ref{P:AndrewsPartialTheta}, which improved some of our identities, and for many discussions about these identities. The views expressed in this article are those of the authors and do not reflect the official policy or position of the U.S. Naval Academy, Department of the Navy, the Department of Defense, or the U.S. Government.
	
	\section{Preliminaries} \label{S: Preliminaries}
	
	In this section, we recall various identities we use in this paper. Andrews' book \cite{AndrewsBook} or Fine's book \cite{Fine} are excellent introductions to $q$-series transformations of these types.
	
	\subsection{The Heine transformation}
	
	We begin by stating the Heine transformation, which is one of the most foundational transformation formulas in the theory of $q$-hypergeometric series.
	
	\begin{lem}[{\cite[Corollary 2.3]{AndrewsBook}}] \label{L: Heine}
		For $|q|,|t|<1$ and $0 < |b| < 1$, we have
		\begin{align*}
			\sum_{n\geq0} \frac{\lp a, b\rp _n}{\lp q, c\rp _n} t^n = \frac{\lp b, at\rp _\infty}{\lp c, t\rp _\infty} \sum_{n\geq0} \frac{\lp \frac{c}{b}, t\rp _n}{\lp q, at\rp _n} b^n.
		\end{align*}
	\end{lem}
	
	We also use the iterated Heine's transformation\footnote{The analytic conditions can be easily derived from those in Lemma \ref{L: Heine}.} (see \cite{GR}, (III.2)). 
	\begin{prop}\label{IteratedHeine} For $|q|, |t| <1$ and $ 0 < |c| < |b| < 1$, we have
		\begin{align*}
			\sum_{n\ge0} \frac{(a,b)_n}{(q,c)_n} t^n = \dfrac{\left(\frac{c}{b},bt\right)_\infty}{(c,t)_\infty} \sum_{n\ge0} \frac{\lp \frac{abt}{c},b\rp_n}{(q,bt)_n} \lp \dfrac{c}{b} \rp^n .
		\end{align*}
	\end{prop}
	
	Finally, we require the following identity, which is derived by Andrews, Subbarao, and Vidyasagar as a consequence of Lemma \ref{L: Heine}. 
	
	\begin{prop}[{\cite[equation (4.1)]{AndrewsSubbaraoVidyasagar}}] \label{P: Andrews-Subbarao-Vidyasagar}
		We have\footnote{Although we do not specify the range in which this and the following identities hold, they are straightforwardly derived from the proof in \cite{AndrewsSubbaraoVidyasagar}.}
		\begin{align*}
			\sum_{n \geq 0} \dfrac{(x)_n}{(y)_n} q^n = \dfrac{q(x)_\infty}{y\lp 1 - \frac{xq}{y} \rp (y)_\infty} + \dfrac{1 - \frac{q}{y}}{1 - \frac{xq}{y}}.
		\end{align*}
	\end{prop}
	
	\subsection{Identities from the Lost Notebook}
	
	A final identity which we require emerges from the early work on Ramanujan's Lost Notebook. In particular, we use the following very general transformation formula due to Andrews that is used frequently to prove identities for partial theta functions.
	
	\begin{prop}[{\cite[Theorem 1]{AndrewsPartialTheta}}] \label{P:AndrewsPartialTheta}
		We have 
		\begin{align*}
			&\sum_{n\ge0} \frac{(B, -Abq)_n q^n}{(-aq, -bq)_n}\label{AndrewsPartialTheta}= \frac{-a^{-1} (B) _\infty (-Abq) _\infty}{(-bq) _\infty( -aq) _\infty} \sum_{n\ge0} \frac{\lp A^{-1}\rp _n}{\lp  - \frac{B}{a} \rp _{n+1}}\lp  \frac{Abq}{a} \rp ^{n}\\
			&\hspace{6cm}+ (1+b)  \sum_{n \ge0} \frac{\lp -a^{-1}\rp _{n+1} \lp  - \frac{ABq}{a} \rp _n }{\lp -\frac{B}{a}, \frac{Abq}{a} \rp _{n+1}}\lp -b\rp ^n.\nonumber
		\end{align*}
	\end{prop}
	
	We also require the following formula from Andrews' work, which is a corollary he derived from Proposition \ref{P:AndrewsPartialTheta}.
	
	\begin{prop}[{\cite[equation (3.9)]{AndrewsPartialTheta}}] \label{P:AndrewsCorollary}
		We have
		\begin{multline*}
			\sum_{n \geq 0} \dfrac{\lp -Bq^2, -Aq^2; q^2 \rp_n}{\lp -aq^2;q^2 \rp_n} q^{2n} = \dfrac{-a^{-1} \lp -Bq^2, -Aq^2; q^2 \rp_\infty}{\lp -aq^2;q^2 \rp_\infty} \sum_{n \geq 0} \dfrac{\lp \frac{Aq^2}{a} \rp^n}{\lp \frac{Bq^2}{a}; q^2 \rp_{n+1}} \\ + \sum_{n \geq 0} \dfrac{\lp -a^{-1}; q^2 \rp_{n+1}}{\lp \frac{Bq^2}{a}, \frac{Aq^2}{a}; q^2 \rp_{n+1}} \lp \frac{AB}{a} \rp^n q^{n^2+3n}.
		\end{multline*}
	\end{prop}
	
	\section{Proof of Theorem \ref{T: Standard}} \label{S: Standard}
	
	In this section, we prove all claimed identities for the functions $\overline{F}_{\x\y}^{\z\w}$.
	
	\subsection{Proof of part {\rm(1)}}\label{S:OU-EU}
	We may build up a generating function for $\overline{F}_\eu^\ou$ by summing over the cases where $2n$ is the size of the maximum even part of the overpartition for each $n \in \mathbb{N}$. There are two possibilities: either the overpartition has a marked part of size $2n$ or it only has unmarked parts of size $2n$. Both cases yield the same expression for the generating function, so the construction is a multiple of two. This gives
	\begin{align*}
		\overline{F}_{\text{eu}}^{\text{ou}}(q) &= 2 \frac{\lp -q;q^2\rp _\infty}{\lp q;q^2\rp _\infty} \sum_{n \geq 0 } \frac{\lp -q^2; q^2\rp _{n-1} \lp q;q^2\rp _n}{\lp q^2, -q; q^2\rp _n} q^{2n} \\&=
		 \frac{\lp -q;q^2\rp _\infty} {\lp q;q^2\rp _\infty} \sum_{n \geq 0 } \frac{\lp -1, q;q^2\rp _n}{\lp q^2, -q; q^2\rp _n} q^{2n}.
	\end{align*}
	Using Lemma \ref{L: Heine} with $q \mapsto q^2$, $a = -1$, $b = q$, $c = -q$, and $t = q^2$, we obtain
	\begin{equation*}
		\overline{F}_{\text{eu}}^{\text{ou}}(q) = 2\frac{\lp -q^2; q^2\rp _\infty}{\lp q^2; q^2\rp _\infty} \sum_{n \geq 0} \frac{q^n}{1 + q^{2n}}.
	\end{equation*}
	To finish the proof, we write
	\begin{equation*}
		2\sum_{n \geq 0} \frac{q^n}{1 + q^{2n}} = 1 + 2 \sum_{n \geq 1} \lp  \sum_{\substack{d \mid n \\ d \equiv 1 \pmod{4}}} 1 - \sum_{\substack{d \mid n \\ d \equiv 3 \pmod{4}}} 1 \rp  q^n.
	\end{equation*}
	It is well-known that
	\begin{equation*}
		\sum_{\substack{d \mid n \\ d \equiv 1 \pmod{4}}} 1 - \sum_{\substack{d \mid n \\ d \equiv 3 \pmod{4}}} 1 =  \frac{r_2(n)}{4},
	\end{equation*}
	where $r_2(n)$ is the number of ways to express $n$ as the sum of squares of two integers. Noting that
	\begin{equation*}
		\sum_{n\ge0} r_2(n)q^n = \Theta^2(\tau)
	\end{equation*}
	gives the claim. 
	
	\subsection{Proof of part {\rm(2)}}\label{S:EU-OU}
	The construction of the generating function $\overline{F}_{\text{ou}}^{\text{eu}}(q)$ is analogous to that of $\overline{F}_{\text{eu}}^{\text{ou}}(q)$. Namely, we sum over the instances where the largest part is odd, and we distinguish two cases based on whether there exists an overlined part of that size. As with $\overline{F}_{\text{eu}}^{\text{ou}}(q)$, these cases turn out the same generating function. Moreover, we account separately for the situation where there are no odd parts with the leading term. We obtain the identity
	\begin{equation*}
		\overline{F}_{\text{ou}}^{\text{eu}}(q) = \dfrac{\lp -q^2;q^2 \rp_\infty}{\lp q^2;q^2 \rp_\infty} + \frac{2q}{1-q} \frac{\lp -q^2; q^2\rp _\infty}{\lp q^2; q^2\rp _\infty} \sum_{n \geq 0} \frac{\lp -q, q^2; q^2\rp _{n}}{\lp q^3,-q^2; q^2\rp_{n}} q^{2n}.
	\end{equation*}
	This completes the proof.
	
	\subsection{Proof of part {\rm(3)}}
	
	We again follow Subsection \ref{S:OU-EU}, splitting into two cases where the largest even part is either marked or not and including a leading term to account for the case where there are no even parts, and we obtain
	\begin{align*}
		\overline{F}_{\text{ed}}^{\text{od}}(q) &= \left(-2q; q^2\right)_\infty + 2 q^2 \lp -2q^3; q^2\rp _\infty \sum_{n\ge 0} \frac{\lp -2q^2; q^2\rp _n}{\lp -2q^3; q^2\rp _n} q^{2n}.
	\end{align*}
	Using Proposition \ref{P: Andrews-Subbarao-Vidyasagar} with $q \mapsto q^2$ and then $x = -2q^2$ and $y = -2q^3$, gives the claim. 
	
	\subsection{Proof of part {\rm(4)}}
	
	We again follow Subsection \ref{S:OU-EU}, splitting into two cases depending on whether the largest odd part is either marked or not and including a leading term to account for the case where there are no odd parts, and we obtain
	\begin{equation*}
		\overline{F}_{\text{od}}^{\text{ed}}(q) = \left(-2q^2; q^2\right) _\infty + 2 q \lp -2q^2; q^2\rp _\infty \sum_{n \ge 0} \frac{\lp -2q; q^2\rp _n}{\lp -2q^2; q^2\rp _n} q^{2n}.
	\end{equation*}
	Using Proposition \ref{P: Andrews-Subbarao-Vidyasagar} with $q \mapsto q^2$, and then $x = -2q$ and $y = - 2q^2$, we obtain the claim.
	
	\subsection{Proof of part {\rm(5)}}\label{S:OD-EU}
	
	We again follow Subsection \ref{S:OU-EU}, splitting into two cases depending on whether the largest even part is either marked or not, and obtain 
	\begin{align*}
		\overline{F}_\eu^\od(q) = 2 \sum_{n \geq 0} \dfrac{\lp -q^2;q^2 \rp_{n-1} \lp -2q^{2n+1};q^2 \rp_\infty}{\lp q^2;q^2 \rp_n} q^{2n} = \lp -2q;q^2 \rp_\infty \sum_{n \geq 0} \dfrac{\lp -1;q^2 \rp_n}{\lp -2q,q^2;q^2 \rp_n} q^{2n}.
	\end{align*}
	Letting $q \mapsto q^2$ and then $t = -\frac{2q}{a}$ in Lemma \ref{L: Heine}, we obtain
	\begin{align*}
		\sum_{n \geq 0} \dfrac{\lp a, b;q^2 \rp_n}{\lp q^2, c; q^2 \rp_n} \left(\frac{-2q}{a}\right)^n = \dfrac{\lp b, - 2q;q^2 \rp_\infty}{\lp c, -\frac{2q}{a}; q^2 \rp_\infty} \sum_{n \geq 0} \dfrac{\lp \frac{c}{b}, -\frac{2q}{a}; q^2 \rp_n}{\lp q^2, -2q;q^2 \rp_n} b^n.
	\end{align*}
	Using 
	\begin{equation}\label{alimit}
		\lim_{a \to \infty} \frac{(a;q)_n}{a^n}=(-1)^n q^{\frac{n(n-1)}{2}},	
	\end{equation}
	we obtain, letting $a \to \infty$, that
	\begin{align*}
		\sum_{n \geq 0} \dfrac{2^n \lp b;q^2 \rp_n q^{n^2}}{\lp q^2, c;q^2 \rp_n} = \dfrac{\lp b, -2q;q^2 \rp_\infty}{\lp c;q^2 \rp_\infty} \sum_{n \geq 0} \dfrac{\lp \frac{c}{b}; q^2 \rp_n}{\lp q^2, -2q;q^2 \rp_n} b^n.
	\end{align*}
	Finally, letting $b = q^2$ and $c = -q^2$ it follows that
	\begin{align*}
		\sum_{n \geq 0} \dfrac{2^n q^{n^2}}{\lp -q^2;q^2 \rp_n} = \dfrac{\lp q^2, -2q;q^2 \rp_\infty}{\lp -q^2;q^2 \rp_\infty} \sum_{n \geq 0} \dfrac{\lp -1;q^2 \rp_n}{\lp q^2,-2q;q^2 \rp_n} q^{2n}.
	\end{align*}
	Rearranging and applying the previous formula for $\overline{F}_\eu^\od(q)$, we conclude the claimed formula for $\overline{F}_\eu^\od(q)$.
	
	\subsection{Proof of part {\rm(6)}}
	
	We again follow Subsection \ref{S:OU-EU}, splitting into two cases depending on whether the largest odd part is either marked or not and including a leading term to account for the case where there are no odd parts. We obtain
	\begin{align*}
		\overline{F}_{\text{od}}^{\text{eu}}(q)
		&=\frac{\left(-q^2; q^2\right)_\infty}{\left(q^2; q^2\right)_\infty} + 2 \sum_{n \ge 0} \frac{\lp -2q; q^2\rp _n \lp -q^{2n+2}; q^2\rp _\infty}{\lp q^{2n+2}; q^2\rp _\infty} q^{2n+1} \\
		&= \frac{\left(-q^2; q^2\right)_\infty}{\left(q^2; q^2\right)_\infty} + 2q \frac{\lp -q^2; q^2\rp _\infty}{\lp q^2; q^2\rp _\infty} \sum_{n \ge 0} \frac{\lp -2q,q^2; q^2\rp _n}{\lp -q^2; q^2\rp _n} q^{2n}.
	\end{align*}
	Applying Proposition \ref{P:AndrewsCorollary} with $a=1$, $A = -1$, and $B = \frac{2}{q}$, gives part (6).
	
	\subsection{Proof of part {\rm(7)}}
	
	We again follow Subsection \ref{S:OU-EU}, splitting into two cases depending on whether the largest even part is either marked or not and including a leading term to account for the case where there are no even parts. We obtain
	\begin{align*}
		\overline{F}_{\text{ed}}^{\text{ou}}(q) &= \frac{\left(-q^2; q^2\right)_\infty}{\left(q^2; q^2\right)_\infty} + 2 \sum_{n \ge 0} \frac{\lp -2q^2; q^2\rp _n \lp -q^{2n+3}; q^2\rp _\infty}{\lp q^{2n+3}; q^2\rp _\infty} q^{2n+2},\\
		&= \frac{\left(-q^2; q^2\right)_\infty}{\left(q^2; q^2\right)_\infty} + 2q^2 \frac{\lp -q^3; q^2\rp _\infty}{\lp q^3; q^2\rp _\infty} \sum_{n \ge 0} \frac{\lp -2q^2, q^3; q^2\rp _n}{\lp -q^3; q^2\rp _n} q^{2n}.
	\end{align*}
	
	Using Proposition \ref{P:AndrewsCorollary} with $B=2$, $A = -q$, and $a=q$, we obtain
	\begin{align*}
		\sum_{n \geq 0} \dfrac{\lp -2q^2, q^3; q^2 \rp_n}{\lp -q^3; q^2 \rp_n} q^{2n} &= - \dfrac{\lp -2q^2, q^3; q^2 \rp_\infty}{q\lp -q^3; q^2 \rp_\infty} \sum_{n \geq 0} \dfrac{(-1)^n q^{2n}}{\lp 2q; q^2 \rp_{n+1}} \\ &\hspace{1cm}+ \lp 1+\dfrac{1}{q} \rp \sum_{n \geq 0} \dfrac{(-1)^n 2^n \lp-q;q^2 \rp_n  q^{n^2+3n}}{\lp 2q, -q^2; q^2 \rp_{n+1}}.
    \end{align*}
    Plugging this into the formula for $\overline{F}_{\text{ed}}^{\text{ou}}(q)$ then gives the claim.
	
	\subsection{Proof of part {\rm(8)}} We again follow Subsection \ref{S:OU-EU}, splitting into two cases depending on whether the largest odd part is either marked or not and including a leading term to account for the case where there are no odd parts. This quickly yields the claimed result, namely
	\begin{align*}
		\overline{F}_{\ou}^{\ed}(q) &= \lp -2q^2;q^2 \rp_\infty + \sum_{n \geq 0} \dfrac{\lp -q;q^2 \rp_n \lp -2q^{2n+2};q^2 \rp_\infty}{\lp q;q^2 \rp_{n+1}} q^{2n+1}\\ &= \lp -2q^2;q^2 \rp_\infty \lp 1 + \dfrac{2q}{1-q} \sum_{n \geq 0} \dfrac{\lp  -q;q^2 \rp _n}{\lp  q^3, -2q^2;q^2 \rp _n} q^{2n} \rp.
	\end{align*} 
	
	\section{Proof of Theorem \ref{T: Modified}} \label{S: Modified}
	
	In this section, we prove Theorem \ref{T: Modified}.
	
	\subsection{Proof of part {\rm(1)}} \label{S: OD-EU mod}
	We modify our argument from Subsection \ref{S:OU-EU} and obtain
	\begin{align*}
		\underline{F}_{\eu}^{\od}(q) = 2 \lp -q;q^2 \rp_\infty \sum_{n \geq 0} \dfrac{\lp -q^2;q^2 \rp_{n-1}}{\lp q^2, -q;q^2 \rp_n} q^{2n} = \lp -q;q^2 \rp_\infty \sum_{n \geq 0} \dfrac{\lp -1;q^2 \rp_n}{\lp q^2, -q;q^2 \rp_n} q^{2n}.
	\end{align*}
	We again use Lemma \ref{L: Heine}. Letting $q \mapsto q^2$ and then $t = - \frac{q}{a}$ we obtain
	\begin{equation}\label{E:SumQABCD}
		\sum_{n \geq 0} \dfrac{\lp a, b;q^2 \rp_n}{\lp q^2, c;q^2 \rp_n} \left(\frac{-q}{a}\right)^n = \dfrac{\lp b, -q;q^2 \rp_\infty}{\lp c, -\frac{q}{a};q^2 \rp_\infty} \sum_{n \geq 0} \dfrac{\lp \frac{c}{b}, -\frac{q}{a};q^2 \rp_n}{\lp q^2, -q;q^2 \rp_n} b^n.
	\end{equation}
	Using \eqref{alimit}, we obtain
	\begin{align*}
		\sum_{n \geq 0} \dfrac{\lp b;q^2 \rp_n q^{n^2}}{\lp q^2, c;q^2 \rp_n} = \dfrac{\lp b, -q;q^2 \rp_\infty}{\lp c;q^2 \rp_\infty} \sum_{n \geq 0} \dfrac{\lp \frac{c}{b}; q^2 \rp_n}{\lp q^2, -q;q^2 \rp_n} b^n.
	\end{align*}
	Setting $b = q^2$ and $c = -q^2$, we conclude 
	\begin{align*}
		\phi(q) = \sum_{n \geq 0} \dfrac{q^{n^2}}{\lp -q^2;q^2 \rp_n} = \dfrac{\lp q^2, -q;q^2 \rp_\infty}{\lp -q^2;q^2 \rp_\infty} \sum_{n \geq 0} \dfrac{\lp -1;q^2 \rp_n}{\lp q^2, -q;q^2 \rp_n} q^{2n} = \dfrac{\lp q^2;q^2 \rp_\infty}{\lp -q^2;q^2 \rp_\infty} \underline{F}_{\eu}^\od(q),
	\end{align*}
	from which the claim follows.
	
	\subsection{Proof of part {\rm(2)}}
	
	As above, we have
	\begin{align*}
		\underline{F}_{\od}^{\eu}(q) = \dfrac{\lp -q^2;q^2 \rp_\infty}{\lp q^2;q^2 \rp_\infty} \lp 1 + 2 q \sum_{n \geq 0} \dfrac{\lp -q;q^2 \rp_n \lp q^2;q^2 \rp_n}{\lp -q^2;q^2 \rp_n} q^{2n} \rp.
	\end{align*}
	Using Proposition \ref{P:AndrewsCorollary} with $B = \frac{1}{q}$, $A = -1$, and $a=1$, we have
	\begin{align*}
		\sum_{n \geq 0} \dfrac{\lp -q, q^2; q^2 \rp_n}{\lp -q^2; q^2 \rp_n} q^{2n} &= - \dfrac{\lp -q, q^2; q^2 \rp_\infty}{\lp -q^2; q^2 \rp_\infty} \sum_{n \geq 0} \dfrac{(-1)^n q^{2n}}{\lp q;q^2 \rp_{n+1}} + 2 \sum_{n \geq 0} \dfrac{(-1)^n q^{n^2 + 2n}}{\lp 1 + q^{2n+2} \rp \lp q;q^2 \rp_{n+1}}.
	\end{align*}
	We next claim that the first summation can be written
	\begin{align*}
		\sum_{n \geq 0} \dfrac{(-1)^n q^{2n}}{\lp q;q^2 \rp_{n+1}} = \sum_{n\ge0}\frac{q^n}{\left(-q^2;q^2\right)_{n+1}}.
	\end{align*}
	Changing variables $q \mapsto q^{\frac 12}$, the claim is equivalent to
	\begin{align*}
		 \frac{1}{1-q^{\frac{1}{2}}} \sum_{n\ge0} \frac{(-1)^n q^n}{\left(q^{\frac{3}{2}};q\right)_n} = \frac{1}{1+q} \sum_{n\ge0} \frac{q^{\frac{n}{2}}}{\left(-q^2;q\right)_n}.
	\end{align*}
	Using Proposition \ref{IteratedHeine} with $a = 0, b = q, c = -q^2$, and $t=q^{\frac 12}$, we obtain the claim. In turn, we complete the proof.
	
	\subsection{Proof of part {\rm(3)}}
	
	We follow the argument in Subsection \ref{S: OD-EU mod}. Namely, we split into two terms based on whether the largest even part is marked or not and include a leading term to account for the case where there are no even parts. This gives
	\begin{align*}
		\underline{F}_{\ed}^{\ou}(q) &= \dfrac{\lp -q;q^2 \rp_\infty}{\lp q;q^2 \rp_\infty} \left( 1 +  2 \sum_{n \geq 1} \dfrac{\lp q;q^2 \rp_n \lp -q^2;q^2 \rp_{n-1}}{\lp -q;q^2 \rp_n} q^{2n} \right) \\
		&= \dfrac{\lp -q;q^2 \rp_\infty}{\lp q;q^2 \rp_\infty} \sum_{n \geq 0} \dfrac{\lp q, -1;q^2 \rp_n}{\lp -q;q^2 \rp_n} q^{2n}.
	\end{align*}
	Using Proposition \ref{P:AndrewsCorollary} with $B = -\frac{1}{q}$, $A=\frac{1}{q^2}$, and $a = \frac{1}{q}$, we have 
	\begin{align*}
		\sum_{n \geq 0} \dfrac{\lp q, -1; q^2 \rp_n}{\lp -q;q^2 \rp_n} q^{2n} &= - \dfrac{q\lp q, -1; q^2 \rp_\infty}{\lp -q;q^2 \rp_\infty} \sum_{n \geq 0} \dfrac{q^{n}}{\lp -q^2;q^2 \rp_{n+1}}\\
		&\hspace{5cm} + \sum_{n \geq 0} \dfrac{(-1)^n \lp -q;q^2 \rp_{n+1} q^{n^2+n}}{\lp -q^2, q; q^2 \rp_{n+1}},
	\end{align*}
        from which the desired result follows.
	
	\subsection{Proof of part {\rm(4)}}
	
	We follow the argument in Subsection \ref{S: OD-EU mod}. Namely, we split into two terms based on whether the largest odd part is marked or not and include a leading term to account for the case where there are no even parts. This directly gives the claim.


\end{document}